\documentclass[final,12pt]{elsarticle}
\usepackage{ragged2e}
\justifying
\usepackage{bm}
\usepackage{geometry}
\usepackage{amsmath}
\usepackage{mathrsfs}
\usepackage{amsthm}
\usepackage{amssymb}
\newtheorem{theorem}{Theorem}[section]
\newtheorem{corollary}{Corollary}[section]
\newtheorem{lemma}{Lemma}[section]

\newenvironment{proof of Theorem 1.1}{\noindent{\textbf{Proof of theorem 1.1.}}\ }{\hfill $\square$\par}

\numberwithin{equation}{section}
\newenvironment{proof1}{\noindent{\textbf{Proof.}}\ }{\hfill $\square$\par}

\usepackage{indentfirst}
\begin{document}
\begin{frontmatter}
  \title{On the maximal $\alpha$-spectral radius of graphs with given matching number\,\tnoteref{titlenote}}
  \tnotetext[titlenote]{This work was supported by the National Nature Science Foundation
   of China (Nos. 11871040).}
  
  \author{Xiying Yuan}
  \ead{xiyingyuan@shu.edu.cn}
 
  \author{Zhenan Shao\corref{correspondingauthor}}
  \cortext[correspondingauthor]{Corresponding author. Email: {\tt shaozhenan@gmail.com} (Zhenan Shao).}
  %\ead{leliu@usst.edu.cn}
  
  \address{Department of Mathematics, Shanghai University, Shanghai 200444, P.R. China}
      
  \begin{abstract}
\indent Let $\mathscr{G}_{n,\beta}$ be the set of graphs of order $n$ with given matching number $\beta$. Let $D(G)$ be the diagonal matrix of the degrees of the graph $G$ and $A(G)$ be the adjacency matrix of the graph $G$. The largest eigenvalue of the nonnegative matrix $A_{\alpha}(G)=\alpha D(G)+A(G)$ is called the $\alpha$-spectral radius of $G$. The graphs with maximal $\alpha$-spectral radius in $\mathscr{G}_{n,\beta}$ are completely characterized in this paper. In this way we provide a general framework to attack the problem of extremal spectral radius in $\mathscr{G}_{n,\beta}$. More precisely, we generalize the known results on the maximal adjacency spectral radius in $\mathscr{G}_{n,\beta}$ and the signless Laplacian spectral radius.

  \end{abstract}
  
  \begin{keyword}
  Adjacency matrix\sep
  Signless Laplacian matrix\sep
   $\alpha$-spectral radius\sep
   Matching number\sep
   Extremal graph
   \MSC[2010]
   05C50
  \end{keyword}
 \end{frontmatter}
%chapter 1
\section{Introduction}
\indent All graphs considered here are simple and undirected. Let $G$ be a graph with vertex set $V(G)$ and edge set $E(G)$. Two distinct edges in a graph $G$ are independent if they are not incident with a common vertex. A set of pairwise independent edges in $G$ is called a \emph{matching} of $G$. The \emph{matching number} $\beta(G)$ (or just $\beta$ for short) of $G$ of order $n$ is the cardinality of maximum matching of $G$, usually, $\beta(G)\leq\frac{n}{2}$. If $\beta(G)=\frac{n}{2}$, we say $G$ has a \emph{perfect matching}. For any $S\subseteq V(G)$, the subgraph induced by $S$ is denoted by $G[S]$. The complete graph of order $n$ is denoted by $K_{n}$, and its complement is denoted by $\overline{K_{n}}$. Let $G_{1}=(V_{1},E_{1})$ and $G_{2}=(V_{2},E_{2})$ be two graphs. The \emph{union} of $G_{1}$ and $G_{2}$ is $G_{1}\cup G_{2}=(V_{1}\cup V_{2}, E_{1}\cup E_{2})$. The \emph{join} of $G_{1}$ and $G_{2}$ is obtained from $G_{1}\cup G_{2}$ by joining edges from each vertex of $G_{1}$ to each vertex of $G_{2}$. The \emph{components} of a graph $G$ are its maximal connected subgraphs. Components with odd (even) order are called the \emph{odd (even) components}. \\
\indent The adjacency matrix of $G$ is $A(G)=(a_{ij})$, where $a_{ij}$=1 if $v_{i}$ is adjacent to $v_{j}$ and 0 otherwise. The signless Laplacian matrix of $G$ is $Q(G)=D(G)+A(G)$ where $D(G)$=diag$(d_{v}:v\in V(G))$. We call the largest eigenvalue of $A(G)$ the \emph{adjacency spectral radius} of $G$ and the largest eigenvalue of $Q(G)$ the \emph{signless Laplacian spectral radius} of $G$. Nikiforov associated the nonnegative matrix $\alpha D(G)+(1-\alpha)A(G)$ with the graph $G$ to carry on the research of eigenvalues and graphs \cite{Niki}. Recently, Liu et al.\cite{liu} modified it to another alternative nonnegative matrix $\alpha D(G)+A(G)$. Write $A_{\alpha}(G)=\alpha D(G)+A(G)$ where $\alpha\geq 0$, and the largest eigenvalue of $A_{\alpha}(G)$ is called the \emph{$\alpha$-spectral radius} of $G$ in this paper, denoted by $\rho_{\alpha}(G)$. \\
\indent The graphs with maximal adjacency spectral radius with given size were characterized in \cite{row}. Zhai et al. characterized the graphs with maximal signless Laplacian spectral radius with given size, clique number or chromatic number in \cite{zhai}. Li et al. generalized these conclusions to the $\alpha$-spectral radius \cite{dan}. The graphs with minimal adjacency spectral radius with given diameter were characterized in \cite{dam}. Fan et al. characterized the graphs with minimal signless Laplacian spectral radius with given number of pendent vertices \cite{fan}. Chen et al. proved that a graph of order sufficiently large $n$ with the minimal degree is Hamilton-connected except for two classes of graphs with given lower bound on the number of edges \cite{chen}. In this paper, we focus on the graphs with given matching number. Let $\mathscr{G}_{n,\beta}$ be the set of graphs of order $n$ with given matching number $\beta$. We characterize the graphs with maximal $\alpha$-spectral radius in $\mathscr{G}_{n,\beta}$ and the main results are as follows.
\begin{theorem}
Let $G$ be any graph in $\mathscr{G}_{n,\beta}$.\\
(1)If $n=2\beta$ or $n=2\beta+1$, then $\rho_{\alpha}(G)\leq(\alpha+1)(n-1)$ with equality if and only if $G\cong K_{n}$.\\
(2)If $2\beta+2\leq n<\frac{(2\alpha+3)\beta+\alpha+2}{\alpha+1}$, then $\rho_{\alpha}(G)\leq2(\alpha+1)\beta$ with equality if and only if $G\cong K_{2\beta+1}\cup\overline{K_{n-2\beta-1}}$.\\
(3)If $n=\frac{(2\alpha+3)\beta+\alpha+2}{\alpha+1}$, then $\rho_{\alpha}(G)\leq2(\alpha+1)\beta$ with equality if and only if $G\cong K_{\beta}\vee\overline{K_{n-\beta}}$ or $G\cong K_{2\beta+1}\cup\overline{K_{n-2\beta-1}}$.\\
(4)If $n>\frac{(2\alpha+3)\beta+\alpha+2}{\alpha+1}$, then 
\begin{align*}
\rho_{\alpha}(G)\leq&\frac{\sqrt{\big[\alpha n+(\alpha+1)\beta-(\alpha+1)\big]^{2}-4\big[(\alpha^{2}-1)\beta n+(\alpha+1)\beta^{2}-\alpha(\alpha+1)\beta\big]}}{2}\\
&+\frac{\alpha n+(\alpha+1)\beta-(\alpha+1)}{2}
\end{align*}
with equality if and only if $G\cong K_{\beta}\vee\overline{K_{n-\beta}}$.
\end{theorem}
\indent Feng et al. considered the graphs with maximal adjacency spectral radius in $\mathscr{G}_{n,\beta}$ (see \cite{feng}) and Yu characterized the graphs with maximal signless Laplacian spectral radius in $\mathscr{G}_{n,\beta}$ (see \cite{yu}). Obviously, $\rho_{0}(G)$ is the adjacency spectral of graph $G$ and $\rho_{1}(G)$ is the signless Laplacian spectral radius of $G$. In fact, we provide a general framework to attack the problem of maximal spectral radius in $\mathscr{G}_{n,\beta}$. More precisely, the following two conclusions are obtained.
\begin{corollary}(\cite[Theorem 1.1]{feng})
Let $G$ be any graph in $\mathscr{G}_{n,\beta}$.\\
(1)If $n=2\beta$ or $n=2\beta+1$, then $\rho_{0}(G)\leq n-1$ with equality if and only if $G\cong K_{n}$.\\
(2)If $2\beta+2\leq n<3\beta+2$, then $\rho_{0}(G)\leq2\beta$ with equality if and only if $G\cong K_{2\beta+1}\cup\overline{K_{n-2\beta-1}}$.\\
(3)If $n=3\beta+2$, then $\rho_{0}(G)\leq2\beta$ with equality if and only if $G\cong K_{\beta}\vee\overline{K_{n-\beta}}$ or $G\cong K_{2\beta+1}\cup\overline{K_{n-2\beta-1}}$.\\
(4)If $n>3\beta+2$, then $\rho_{0}(G)\leq \frac{1}{2}\big(\beta-1+\sqrt{(\beta-1)^{2}+4\beta(n-\beta)}\big)$ with equality if and only if $G\cong K_{\beta}\vee\overline{K_{n-\beta}}$.
\end{corollary}
\begin{corollary}(\cite[Theorem 2.5]{yu})
Let $G$ be any graph in $\mathscr{G}_{n,\beta}$.\\
(1)If $n=2\beta$ or $n=2\beta+1$, then $\rho_{1}(G)\leq 2(n-1)$ with equality if and only if $G\cong K_{n}$.\\
(2)If $2\beta+2\leq n<\frac{5\beta+3}{2}$, then $\rho_{1}(G)\leq 4\beta$ with equality if and only if $G\cong K_{2\beta+1}\cup\overline{K_{n-2\beta-1}}$.\\
(3)If $n=\frac{5\beta+3}{2}$, then $\rho_{1}(G)\leq 4\beta$ with equality if and only if $G\cong K_{\beta}\vee\overline{K_{n-\beta}}$ or $G\cong K_{2\beta+1}\cup\overline{K_{n-2\beta-1}}$.\\
(4)If $n>\frac{5\beta+3}{2}$, then $\rho_{1}(G)\leq \frac{1}{2}\big(n-2+2\beta$+$\sqrt{(n-2+2\beta)^{2}-8\beta^{2}+8\beta}\big)$ with equality if and only if $G\cong K_{\beta}\vee\overline{K_{n-\beta}}$.
\end{corollary}
%chapter 2
\section{Auxiliary Results}
{In order to prove the main results, we present several technical lemmas as follows. The first one is the well-known Perron-Frobenius Theorem for irreducible nonnegative matrices. When $G$ is connected, $A_{\alpha}(G)$ is an irreducible nonnegative matrix, and then Perron-Frobenius Theorem is available for $A_{\alpha}(G)$ of a connected graph $G$. Let $\rho(A)$ be the spectral radius of a real symmetric matrix $A$, i.e. $\rho(A)=\max\{|\lambda|:\lambda\ is\ an\ eigenvalue\ of\ A\}$. We write $A\geq B$ if $A-B\geq 0$ (every entry of $A-B$ is not less than 0). 
%lemma2.1
\begin{lemma}(\cite{horn} Perron-Frobenius Theorem)\label{fb}
If A is an irreducible nonnegative matrix of order $n$ with $n\geq2$, then the following statements hold.\\
(1) $A$ has a positive eigenvector corresponding to $\rho(A)$.\\
(2) If $A\geq B$ and $A\neq B$, then $\rho(A)>\rho(B)$.
\end{lemma}
%lemma2.2
\begin{lemma}(\cite{bondy} Tutte-Berge Formula)\label{2.2}
For any graph $G$ of order $n$,
\begin{align*}
\beta(G)=\frac{1}{2}\min\bigg\{n-\big(o(G-S)-|S|\big):S\subset V(G)\bigg\},
\end{align*}
where $o(G-S)$ denotes the number of odd components of $G-S$.
\end{lemma}
\indent For a graph $G$ of order $n$ with given matching number $\beta(G)$, there is a subset $S$ satisfying the Tutte-Berge Formula by Lemma \ref{2.2}. We denote $|S|$ by $s(G)$ (or simply by $s$), and denote by $q(G)=n+s(G)-2\beta(G)$ (or simply by $q$). In virtue of Tutte-Berge Formula, we may give a characterization of the extremal graph $G$, which achieves the maximal $\alpha$-spectral radius in $\mathscr{G}_{n,\beta}$.
%lemma2.3
\begin{lemma}\label{2.3}
Let $G$ be a graph with the maximal $\alpha$-spectral radius in $\mathscr{G}_{n,\beta}$. The following conclusions hold.\\
(1) If $s=0$, $q=0$, then $G\cong K_{2\beta}$.\\
(2) If $s=0$, $q\geq 1$, then $G\cong K_{2\beta+1}\cup\overline{K_{n-2\beta-1}}$.\\
(3) If $s\geq 1$, then $G\cong K_{s}\vee(\mathop\cup\limits^{q}_{i=1}K_{n_{i}})$ where $\sum\limits_{i=1}^{q}n_{i}=n-s$.
\end{lemma}
\begin{proof1}
If $s=0$ and $q=0$, and noting $G$ is the graph with maximal $\alpha$-spectral radius, we have $G\cong K_{2\beta}$. If $s=0$ and $q\geq1$, then $G-S$ has $n-2\beta$ odd components and $m$ $(m\geq 0)$ even components. Since each odd component consists of at least one vertex and each even component (if it exists) consists of at least two vertices, there are at most $2\beta+1$ vertices in any odd component and at most $2\beta$ vertices in any even component. Since $G$ is the graph with maximal $\alpha$-spectral radius,  we deduce $G=K_{2\beta+1}\cup\overline{K_{n-2\beta-1}}$. \\
\indent Now we suppose $s\geq1$. Then $q\geq 1$ holds. Otherwise we have $n=2\beta-s$, which contradicts with $n\geq 2\beta$. Let $G_{1}, G_{2},\cdots, G_{q}$ be the odd components in $G-S$ with $|V(G_{i})|=n_{i}\geq1$ for $i=1,\cdots,q$. Clearly, $n\geq s+q=n+2s-2\beta$. Thus $\beta\geq s$. \\
\indent If $G-S$ has no even components, then $\mathop{\cup}\limits_{i=1}^{q}V(G_{i})=V(G)\backslash S$. The graph with the maximal $\alpha$-spectral radius in $\mathscr{G}_{n,\beta}$ is $G^{\prime}\cong K_{s}\vee(\mathop\cup\limits^{q}_{i=1}K_{n_{i}})$. We have $G=G^{\prime}$. Otherwise $G$ is a proper subgraph of the connected graph $G^{\prime}$. Then by Lemma \ref{fb} (2) we have $\rho_{\alpha}(G^{\prime})>\rho(G)$ while $G^{\prime}\in\mathscr{G}_{n,\beta}$.\\
\indent If $G-S$ contains even components, let $C$ be the union of these even components where $|C|=p$. Let $G^{\prime}$ be the graph obtained from $G$ by adding edges to make $G[V(G_{q})\cup V(C)]$ to be a clique. Then $G^{\prime}$ has no even components. According to the conclusion of the above proof, the graph with maximal $\alpha$-spectral is $G^{\prime\prime}\cong K_{s}\vee(\mathop\cup\limits^{q-1}_{i=1}K_{n_{i}}\cup K_{n_{q}+p})$. Similarly, we have $G=G^{\prime\prime}$.
\end{proof1}
\indent We may further prove that when $s\geq 1$ for the extremal case, some $n_{i}=2\beta-2s+1$, and $n_{j}=1$, $j\neq i$.
%lemma2.4
\begin{lemma}\label{2.4}
Let $G$ be a graph in $\mathscr{G}_{n,\beta}$ with the maximal $\alpha$-spectral radius where $n\geq2\beta+1$. Then
\begin{center}
$G=K_{s}\vee(K_{2\beta-2s+1}\cup\overline{K_{q-1}})$.
\end{center}
\end{lemma}
\begin{proof1}
Let $G_{1}, G_{2},\cdots, G_{q}$ be the odd components in $G-S$ with $|V(G_{i})|=n_{i}\geq1$ for $i=1,\cdots,q$. Without loss of generality, we assume that $n_{1}\leq n_{2}\leq\cdots\leq n_{q}$.\\
\indent When $s\geq 1$, by Lemma \ref{2.3}, we know the graph $G$ with maximal $\alpha$-spectral radius is $G=K_{s}\vee(\mathop\cup\limits^{q}_{i=1}K_{n_{i}})$. Then we only need to prove that $n_{1}=n_{2}=\cdots=n_{q-1}=1$. For convenience, we usually use $\rho_{\alpha}$ to denote $\rho_{\alpha}(G)$. By Lemma \ref{fb} (1), there is a positive eigenvector $\bm x$ associated with $\rho_{\alpha}$. Since $\rho_{\alpha}>0>-1$, the vertices in each odd components admit the same coordinate in $\bm x$ (denoted by $x_{i}$), where $1\leq i\leq q$. The vertices in $S$ admit the same coordinate in $\bm x$ (denoted by $p$). Then $A_{\alpha}(G)\bm x=\rho_{\alpha}\bm x$ implies
\begin{align}
\left\{
\begin{aligned}
&\big[\rho_{\alpha}-(\alpha+1)(n_{1}-1)-\alpha s\big]x_{1}-sp=0,\\
&\big[\rho_{\alpha}-(\alpha+1)(n_{2}-1)-\alpha s\big]x_{2}-sp=0,\\
&\quad\quad\quad\quad\quad\qquad\cdots\\
&\big[\rho_{\alpha}-(\alpha+1)(n_{q}-1)-\alpha s\big]x_{q}-sp=0,\\
&\sum\limits_{i=1}^{q}n_{i}x_{i}-(\rho_{\alpha}-\alpha n-s+\alpha+1)p=0.\label{system}
\end{aligned}
\right.
\end{align}
\indent Let $M$ be the coefficient matrix of system (\ref{system}). Since $(x_{1},x_{2},\cdots,x_{q},p)^{t}$ is a non-zero solution of system (\ref{system}), we have $|M|=0$. By solving $|M|=0$, we get the following relation\\
\begin{align}
|M|=&-\prod\limits_{i=1}^{q}\big[\rho_{\alpha}-(\alpha+1)(n_{i}-1)-\alpha s\big]\notag\\
&\times\bigg[\rho_{\alpha}-\alpha n-s+\alpha+1-\sum\limits_{i=1}^{q}\frac{n_{i}s}{\rho_{\alpha}-(\alpha+1)(n_{i}-1)-\alpha s}\bigg]=0.\label{hanglieshi}
\end{align}
Then
\begin{align*}
\rho_{\alpha}-\alpha n-s+\alpha+1-\sum\limits_{i=1}^{q}\frac{n_{i}s}{\rho_{\alpha}-(\alpha+1)(n_{i}-1)-\alpha s}=0.
\end{align*}
\indent If $n_{q-1}\geq3$, set 
\begin{align*}
G_{0}=K_{s}\vee(K_{n_{q}+2}\cup K_{n_{q-1}-2}\mathop\cup\limits_{i=1}^{q-2}K_{n_{i}}).
\end{align*}
We claim $\rho_{\alpha}(G)<\rho_{\alpha}(G_{0})$. Consider the function
\begin{align*}
f(\delta,\lambda)=&\frac{\lambda-\alpha n-s+\alpha+1}{s}-\sum\limits_{i=1}^{q-2}\frac{n_{i}}{\lambda-(\alpha+1)(n_{i}-1)-\alpha s}\\
&-\frac{n_{q-1}-\delta}{\lambda-(\alpha+1)(n_{q-1}-\delta-1)-\alpha s}-\frac{n_{q}+\delta}{\lambda-(\alpha+1)(n_{q}+\delta-1)-\alpha s}.
\end{align*}
where $\lambda\geq(\alpha+1)(n_{q}+s-1)$, $0\leq\delta\leq2$. Clearly, $f(0,\rho_{\alpha}(G))=0$.\\
\indent Taking derivative with respect to $\delta$, we have
\begin{align*}
\frac{df(\delta,\lambda)}{d\delta}&=\frac{\lambda-\alpha s+(\alpha+1)}{\big[\lambda-(\alpha+1)(n_{q-1}-\delta-1)-\alpha s\big]^{2}}-\frac{\lambda-\alpha s+(\alpha+1)}{\big[\lambda-(\alpha+1)(n_{q}+\delta-1)-\alpha s\big]^{2}}\\
&=\frac{\big[\lambda-\alpha s+(\alpha+1)\big]\big[2\lambda-(\alpha+1)(n_{q}+n_{q-1}-2)-2\alpha s\big](\alpha+1)(n_{q-1}-n_{q}-2\delta)}{\big[\lambda-(\alpha+1)(n_{q-1}-\delta-1)-\alpha s\big]^{2}\big[\lambda-(\alpha+1)(n_{q}+\delta-1)-\alpha s\big]^{2}}\\
&<0.
\end{align*}
Hence, $f(\delta, \lambda)$ strictly decreases with respect to $\delta$ for $\lambda\geq(\alpha+1)(n_{q}+s-1)$. Note that $\rho_{\alpha}(G)\geq(\alpha+1)(n_{q}+s-1)$ and $\rho_{\alpha}(G_{0})\geq(\alpha+1)(n_{q}+s-1)$. So we have $f(2,\rho_{\alpha}(G))<0=f(0,\rho_{\alpha}(G))$.\\
\indent For $f(2,\lambda)$ where $\lambda\geq(\alpha+1)(n_{q}+s-1)$, taking derivative with respect to $\lambda$, we have
\begin{align*}
\frac{df(2,\lambda)}{d\lambda}=&\frac{1}{s}+\sum\limits_{i=1}^{q-2}\frac{n_{i}}{\big[\lambda-(\alpha+1)(n_{i}-1)-\alpha s\big]^{2}}\\
&+\frac{n_{q-1}-2}{\big[\lambda-(\alpha+1)(n_{q-1}-3)-\alpha s\big]^{2}}+\frac{n_{q}+2}{\big[\lambda-(\alpha+1)(n_{q}+1)-\alpha s\big]^{2}}>0.
\end{align*}
Hence $f(2,\lambda)$ increases with respect to $\lambda$ for $\lambda\geq(\alpha+1)(n_{q}+s-1)$.\\
\indent Since $\rho_{\alpha}(G_{0})$ is the largest root of equation $f(2,\lambda)=0$, so $f(2,\rho_{\alpha}(G))<0=f(2,\rho_{\alpha}(G_{0}))$. Thus $\rho_{\alpha}(G)<\rho_{\alpha}(G_{0})$. Moreover, $\beta(G)=\beta(G_{0})$. Using the above arguments repeatedly, we finally get $G=K_{s}\vee(K_{2\beta-2s+1}\cup\overline{K_{q-1}})$.
\end{proof1}
\indent Next, we will calculate the $\alpha$- spectral radius of $K_{\beta}\vee\overline{K_{n-\beta}}$. 
%lemma2.5
\begin{lemma}
For $n>\beta\geq1$,we have
\begin{align}
\rho_{\alpha}(K_{\beta}\vee\overline{K_{n-\beta}})=&\frac{\sqrt{\big[\alpha n+(\alpha+1)\beta-(\alpha+1)\big]^{2}-4\big[(\alpha^{2}-1)\beta n+(\alpha+1)\beta^{2}-\alpha(\alpha+1)\beta\big]}}{2}\notag\\
&+\frac{\alpha n+(\alpha+1)\beta-(\alpha+1)}{2}.\label{2}
\end{align}
\end{lemma}
\begin{proof1}
Note that $K_{\beta}\vee\overline{K_{n-\beta}}=K_{\beta}\vee(K_{1}\cup\overline{K_{n-\beta-1}})$. By Equation (\ref{hanglieshi}), we know that $\rho_{\alpha}(K_{\beta}\vee\overline{K_{n-\beta}})$ satisfies $g(\lambda)=0$, where
\begin{align*}
g(\lambda)=\lambda^{2}-\big[\alpha n+(\alpha+1)\beta-(\alpha+1)\big]\lambda+(\alpha^{2}-1)\beta n+(\alpha+1)\beta^{2}-\alpha(\alpha+1)\beta.
\end{align*}
By a simple calculation, the conclusion follows immediately.
\end{proof1}
%chapter 3
\section{Proof of Theorem 1.1}
\indent The general structure of the graph with maximal $\alpha$-spectra radius in $\mathscr{G}_{n,\beta}$ has been described in section 2 in terms of the parameters $s$ and $q$. Now we will prove Theorem 1.1 by careful computation which explore how $\rho_{\alpha}(G)$ disturbs with the change of $n$ and $s$.\\
\begin{proof of Theorem 1.1}
By Lemma \ref{2.4} and Equation (\ref{hanglieshi}), we know that $\rho_{\alpha}(G)$ satisfies $f(\lambda)=0$, where
\begin{align}\label{f}
f(\lambda)=&(\lambda-\alpha n-s+\alpha+1)(\lambda-\alpha s)\big[\lambda-2(\alpha+1)\beta+(\alpha+2)s\big]\notag\\
&-s(n+s-2\beta-1)\big[\lambda-2(\alpha+1)\beta+(\alpha+2)s\big]-s(2\beta-2s+1)(\lambda-\alpha s).
\end{align}
It is not hard to see that
\begin{align*}
&f(-\infty)<0,\\
&f(\alpha s)=2s(1+\alpha)(\beta-s)(n+s-2\beta-1)\geq 0,\\
&f(2(\alpha+1)\beta-(\alpha+1)s)=s(\alpha+1)(2\alpha\beta-2\alpha s+s)(2\beta-n-s+1)\leq 0,\\
&f(+\infty)>0.
\end{align*}
Hence the three roots of $f(\lambda)=0$ lie in $(-\infty, \alpha s)$, $(\alpha s, 2(\alpha+1)\beta-(\alpha+1)s)$ and $(2(\alpha+1)\beta-(\alpha+1)s,+\infty)$. So $f(\lambda)=0$ has exactly one root not less than $2(\alpha+1)\beta-(\alpha+1)s$.\\
\indent(1) If $n=2\beta$ or $2\beta+1$, it is obvious that $\rho_{\alpha}(G)\leq \rho_{\alpha}(K_{n})$ with equality if and only if $G\cong K_{n}$.\\
\indent(2) If $2\beta+2\leq n<\frac{(2\alpha+3)\beta+\alpha+2}{\alpha+1}$, then we have $\rho_{\alpha}(K_{\beta}\vee\overline{K_{n-\beta}})<\rho_{\alpha}(K_{2\beta+1}\cup\overline{K_{n-2\beta-1}})=2(1+\alpha)\beta$. Moreover,
\begin{align*}
f(2(\alpha+1)\beta)&=\big[2(\alpha+1)\beta-\alpha n-s+\alpha+1\big]\big[2(\alpha+1)\beta-\alpha s\big](\alpha+2)s-(n+s-2\beta-1)s^{2}\\
&\quad\ (\alpha+2)-(2\beta-2s+1)\big[2(\alpha+1)\beta-\alpha s\big]s\\
&\geq\bigg[2(\alpha+1)\beta-\alpha\frac{(2\alpha+3)\beta+\alpha+2}{\alpha+1}-s+\alpha+1\bigg]\big[2(\alpha+1)\beta-\alpha s\big](\alpha+2)s-\\
&\quad\ \bigg[\frac{(2\alpha+3)\beta+\alpha+2}{\alpha+1}+s-2\beta-1\bigg](\alpha+2)s^{2}-(2\beta-2s+1)\big[2(\alpha+1)\beta-\alpha s\big]s\\
&=s\bigg\{\bigg(\frac{\beta}{\alpha+1}+\alpha\beta-\alpha s+\frac{1}{\alpha+1}+\beta\bigg)\big[\alpha(\beta-s)+\alpha\beta+2\beta\big]-\bigg(\frac{\beta}{\alpha+1}+\frac{1}{\alpha+1}\\
&\quad\ +s\bigg)(\alpha s+2s)\bigg\}\\
&\geq s^{2}(\alpha+1)(\alpha+2)(\beta-s)\\
&\geq 0.
\end{align*}
This means that $\rho_{\alpha}(G)\leq 2(\alpha+1)\beta$. If $\rho_{\alpha}(G)=2(\alpha+1)\beta$, then $s=0$ or $s=\beta$. While $\rho_{\alpha}(K_{\beta}\vee\overline{K_{n-\beta}})< 2(\alpha+1)\beta$ when $2\beta+2\leq n<\frac{(2\alpha+3)\beta+\alpha+2}{\alpha+1}$, we deduce $s=0$. By Lemma \ref{2.4}, we have $G\cong K_{2\beta+1}\cup\overline{K_{n-2\beta-1}}$.\\
\indent (3) If $n=\frac{(2\alpha+3)\beta+\alpha+2}{\alpha+1}$, we have $f(2(\alpha+1)\beta)\geq s^{2}(\alpha+1)(\alpha+2)(\beta-s)\geq 0$ which implies $\rho_{\alpha}(G)\leq 2(\alpha+1)\beta$. If $\rho_{\alpha}(G)=2(\alpha+1)\beta$, then $s=0$ or $\beta=s$. By Lemma \ref{2.4}, we have $G\cong K_{\beta}\vee\overline{K_{n-\beta}}$ or $G\cong K_{2\beta+1}\cup\overline{K_{n-2\beta-1}}$.\\
\indent (4) If $n>\frac{(2\alpha+3)\beta+\alpha+2}{\alpha+1}$, then we have $\rho_{\alpha}(K_{\beta}\vee\overline{K_{n-\beta}})>\rho_{\alpha}(K_{2\beta+1}\cup\overline{K_{n-2\beta-1}})=2(\alpha+1)\beta$. Hence, $s\neq 0$. Next, we prove that $\rho_{\alpha}(G)\leq\rho_{\alpha}(K_{\beta}\vee\overline{K_{n-\beta}})$. Note that $\rho_{\alpha}(K_{\beta}\vee\overline{K_{n-\beta}})$ satisfies $g(\lambda)=0$, where
\begin{align*}
g(\lambda)=\lambda^{2}-\big[\alpha n+(\alpha+1)\beta-(\alpha+1)\big]\lambda+(\alpha^{2}-1)\beta n+(\alpha+1)\beta^{2}-\alpha(\alpha+1)\beta.
\end{align*}
We have
\begin{align*}
f(\lambda)=&g(\lambda)\big[\lambda+s-(\alpha+1)\beta\big]+(\alpha+1)(\beta-s)\big[(n+\alpha s-\alpha\beta+s-2\beta-1)\lambda+\alpha^{2}\beta n-a^{2}sn\\
&+2sn-\alpha sn-\beta n-\alpha^{2}\beta+\alpha^{2}s+\alpha\beta^{2}+\alpha\beta s-\alpha\beta-\alpha s^{2}+2\alpha s+b^{2}-4\beta s+2s^{2}-2s\big].
\end{align*}
There are two cases as follows. \\
%case1
\indent \textbf{Case 1.} $n+\alpha s-\alpha\beta+s-2\beta-1\geq 0$.\\
\indent In this case, we have
\begin{align*}
f(\rho_{\alpha}(K_{\beta}\vee\overline{K_{n-\beta}}))&\geq (\alpha+1)(\beta-s)\big[(n+\alpha s-\alpha\beta+s-2\beta-1)2(\alpha+1)\beta+(\alpha^{2}\beta-\alpha^{2}s+\\
&\quad\ 2s-\alpha s-\beta)n-\alpha^{2}\beta+\alpha^{2}s+\alpha\beta^{2}+\alpha\beta s-\alpha\beta-\alpha s^{2}+2\alpha s+b^{2}-4\beta s\\
&\quad\ +2s^{2}-2s\big]\\
&\geq(\alpha+1)(\beta-s)\bigg\{(\alpha s-\alpha\beta+s-2\beta-1)2(1+\alpha)\beta+\big[\beta+2s+(\alpha^{2}\beta-\\
&\quad\ \alpha^{2}s)+(2\alpha\beta-\alpha s)\big]\frac{(2\alpha+3)\beta+\alpha+2}{\alpha+1}-\alpha^{2}\beta+\alpha^{2}s+\alpha\beta^{2}+\alpha\beta s-\alpha\beta\\
&\quad\ -\alpha s^{2}+2\alpha s+b^{2}-4\beta s+2s^{2}-2s\bigg\}\\
&=s(\beta-s)(2\alpha^{2}\beta-\alpha^{2}s+4\alpha\beta+\alpha s+4\beta+2s+2)\geq 0.
\end{align*}
Hence, $f(\rho_{\alpha}(K_{\beta}\vee\overline{K_{n-\beta}}))\geq 0$. This means that $\rho_{\alpha}(G)\leq\rho_{\alpha}(K_{\beta}\vee\overline{K_{n-\beta}})$. If $\rho_{\alpha}(G)=\rho_{\alpha}(K_{\beta}\vee\overline{K_{n-\beta}})$, then $\beta=s$. By Lemma \ref{2.4}, we have $G\cong K_{\beta}\vee\overline{K_{n-\beta}}$.\\
%case2
\indent \textbf{Case 2.} $n+\alpha s-\alpha\beta+s-2\beta-1<0$.\\
\indent Note that $n+\alpha s-\alpha\beta+s-2\beta-1\geq0$ always holds when $0<\alpha\leq\frac{\sqrt{5}-1}{2}$. Then we suppose $\alpha>\frac{\sqrt{5}-1}{2}$. If $\beta=1$, it is obvious that the graph with maximal $\alpha$-spectral is $K_{1}\vee\overline{K_{n-1}}$ when $n>3+\frac{2}{\alpha+1}$. Hence, we suppose $\beta\geq2$. Combining $n<(\alpha+2)\beta-(\alpha+1)s+1$ with $n>\frac{(2\alpha+3)\beta+\alpha+2}{\alpha+1}$, we have $s<\frac{(\alpha^{2}+\alpha-1)\beta}{(1+\alpha)^{2}}-\frac{1}{(1+\alpha)^{2}}<\frac{(\alpha^{2}+\alpha-1)\beta}{(1+\alpha)^{2}}$. Moreover, by Equation (\ref{2}), we deduce $\rho_{\alpha} (K_{\beta}\vee\overline{K_{n-\beta}})\geq \alpha n+\frac{\alpha+2}{\alpha+1}\beta-\frac{\alpha(\alpha+2)}{\alpha+1}$ when $n>\frac{(2\alpha+3)\beta+\alpha+2}{\alpha+1}$ and $s<\frac{(\alpha^{2}+\alpha-1)\beta}{(1+\alpha)^{2}}$.\\
\indent Substituting $\rho_{\alpha}=\alpha n+\frac{\alpha+2}{\alpha+1}\beta-\frac{\alpha(\alpha+2)}{\alpha+1}$ into the Equation (\ref{f}), now we claim that $f(\alpha n+\frac{\alpha+2}{\alpha+1}\beta-\frac{\alpha(\alpha+2)}{\alpha+1})>0$ if $1\leq s\leq\frac{(\alpha^{2}+\alpha-1)\beta}{(1+\alpha)^{2}}$. Taking the second-order partial derivative of the function $f(\alpha n+\frac{\alpha+2}{\alpha+1}\beta-\frac{\alpha(\alpha+2)}{\alpha+1})$ with respect to $s$, we have
\begin{align*}
\frac{\partial ^{2} f(\alpha n+\frac{\alpha+2}{\alpha+1}\beta-\frac{\alpha(\alpha+2)}{\alpha+1})}{\partial s^{2}}=2(\alpha+1)(-3\alpha\beta+3\alpha s+4\beta-2n-6s+2)<0.
\end{align*}
Hence, $f\big(\alpha n+\frac{\alpha+2}{\alpha+1}\beta-\frac{\alpha(\alpha+2)}{\alpha+1}\big)$ reaches the minimum when $s=1$ or $s=\frac{(\alpha^{2}+\alpha-1)\beta}{(1+\alpha)^{2}}$. If $s=1$, then
\begin{align*}
f\bigg(\alpha n+\frac{\alpha+2}{\alpha+1}\beta-\frac{\alpha(\alpha+2)}{\alpha+1}\bigg)=&\frac{1}{(\alpha+1)^{3}}\big[(\alpha^{5}\beta-\alpha^{5}+4\alpha^{4}\beta-3\alpha^{4}+5\alpha^{3}\beta-4\alpha^{3}+2\alpha^{2}\beta-3\alpha^{2}\\
&-\alpha)n^{2}-(2\alpha^{5}\beta^{2}-2\alpha^{5}+8\alpha^{4}\beta^{2}+2\alpha^{4}\beta-5\alpha^{4}+8\alpha^{3}\beta^{2}+\alpha^{3}\beta\\
&-2\alpha^{3}-2\alpha^{2}\beta^{2}-8\alpha^{2}\beta+3\alpha^{2}-4\alpha\beta^{2}-7\alpha\beta+4\alpha+2)n-(-\\
&4\alpha^{5}\beta^{2}+4\alpha^{5}\beta+2\alpha^{4}\beta^{3}-18\alpha^{4}\beta^{2}+10\alpha^{4}\beta+11\alpha^{3}\beta^{3}-25\alpha^{3}\beta^{2}+\\
&4\alpha^{3}\beta+20\alpha^{2}\beta^{3}-6\alpha^{2}\beta^{2}-4\alpha^{2}\beta+2\alpha^{2}+12\alpha\beta^{3}+4\alpha\beta^{2}-5\alpha\beta\\
&+3\alpha-4\beta^{2}-6\beta)\big]
\end{align*}
and
\begin{align*}
\frac{\partial ^{2}f\big(\alpha n+\frac{\alpha+2}{\alpha+1}\beta-\frac{\alpha(\alpha+2)}{\alpha+1}\big)}{\partial ^{2} n}=\frac{2\alpha(\alpha^{2}\beta-\alpha^{2}+2\alpha\beta-\alpha-1)}{\alpha+1}.
\end{align*}
Combining $\beta\geq 2$ with $\alpha>\frac{\sqrt{5}-1}{2}$, we obtain
\begin{align*}
\frac{2\alpha(\alpha^{2}\beta-\alpha^{2}+2\alpha\beta-\alpha-1)}{\alpha+1}\geq\frac{2\alpha(\alpha^{2}+3\alpha-1)}{\alpha+1}>0.
\end{align*}
We deduce that $f\big(\alpha n+\frac{\alpha+2}{\alpha+1}\beta-\frac{\alpha(\alpha+2)}{\alpha+1}\big)$ reaches the minimum $f^{1}_{\min}$ when $n=\frac{(2\alpha+3)\beta+\alpha+2}{\alpha+1}$, where
\begin{align*}
f^{1}_{\min}=(\beta-1)(2\alpha^{2}\beta-\alpha^{2}+4\alpha\beta+\alpha+4\beta+4)>0.
\end{align*}
\indent If $s=\frac{(\alpha^{2}+\alpha-1)\beta}{(1+\alpha)^{2}}$, we have
\begin{align*}
f\bigg(\alpha n+\frac{\alpha+2}{\alpha+1}\beta-\frac{\alpha(\alpha+2)}{\alpha+1}\bigg)=&\frac{1}{(\alpha+1)^{5}}\big[(\alpha^{6}\beta+\alpha^{6}+5\alpha^{5}\beta+4\alpha^{5}+10\alpha^{4}\beta+6\alpha^{4}+10\alpha^{3}\beta+4\alpha^{3}\\
&+5\alpha^{2}\beta+\alpha^{2}+\alpha\beta)n^{2}-(2\alpha^{6}\beta^{2}+5\alpha^{6}\beta+2\alpha^{6}+7\alpha^{5}\beta^{2}+23\alpha^{5}\beta\\
&+10\alpha^{5}+5\alpha^{4}\beta^{2}+40\alpha^{4}\beta+18\alpha^{4}-6\alpha^{3}\beta^{2}+33\alpha^{3}\beta+14\alpha^{3}-7\alpha^{2}\beta^{2}\\
&+13\alpha^{2}\beta+4\alpha^{2}+\alpha\beta^{2}+2\alpha\beta+2\beta^{2})n+(\alpha^{6}\beta^{3}+4\alpha^{6}\beta^{2}+4\alpha^{6}\beta+\\
&\ \alpha^{6}+2\alpha^{5}\beta^{3}+14\alpha^{5}\beta^{2}+21\alpha^{5}\beta+6\alpha^{5}-9\alpha^{4}\beta^{3}+7\alpha^{4}\beta^{2}+39\alpha^{4}\beta\\
&+13\alpha^{4}-39\alpha^{3}\beta^{3}-22\alpha^{3}\beta^{2}+30\alpha^{3}\beta+12\alpha^{3}-56\alpha^{2}\beta^{3}-27\alpha^{2}\beta^{2}\\
&+8\alpha^{2}\beta+4\alpha^{2}-28\alpha\beta^{3}-8\alpha\beta^{2}+2\beta^{3})\big].
\end{align*}
It is not hard to know $f\big(\alpha n+\frac{\alpha+2}{\alpha+1}\beta-\frac{\alpha(\alpha+2)}{\alpha+1}\big)$ reaches the minimum $f^{2}_{\min}$ when $n=\frac{(2\alpha+3)\beta+\alpha+2}{\alpha+1}$, where
\begin{align*}
f^{2}_{\min}=\frac{1}{(\alpha+1)^{5}}\beta^{2}(\alpha+2)(\alpha^{2}+\alpha-1)(\alpha^{3}\beta+7\alpha^{2}\beta+11\alpha\beta+2\alpha+2\beta+2)>0.
\end{align*}
Let $G^{\prime}$ be the graph with maximal $\alpha$-spectral radius with given matching number $\beta$ when $1\leq s\leq\frac{(\alpha^{2}+\alpha-1)\beta}{(1+\alpha)^{2}}$. Since $f(\lambda)$ has exactly one root in $(2(\alpha+1)\beta-(\alpha+1)s,+\infty)$ and $\alpha n+\frac{\alpha+2}{\alpha+1}\beta-\frac{\alpha(\alpha+2)}{\alpha+1}>2(\alpha+1)\beta-(\alpha+1)s$, we have $\rho_{\alpha}(G^{\prime})<\alpha n+\frac{\alpha+2}{\alpha+1}\beta-\frac{\alpha(\alpha+2)}{\alpha+1}\leq\rho_{\alpha}(K_{\beta}\vee\overline{K_{n-\beta}})$.\\
\indent Combining with the above two cases, we deduce that the maximal graph is $K_{\beta}\vee\overline{K_{n-\beta}}$ if $n>\frac{(2\alpha+3)\beta+\alpha+2}{\alpha+1}$. 
\end{proof of Theorem 1.1}
%references

\end{document}